\documentclass [12pt]{article}
\usepackage{graphicx,amssymb,amsfonts,latexsym,amsmath,amsthm,times}
\usepackage{epsfig}
\usepackage{color}
\usepackage[english]{babel}
\usepackage[a4paper]{geometry}
\def\lanbox{\hbox{$\, \vrule height 0.25cm width 0.25cm depth 0.01cm \,$}}

\numberwithin{equation}{section}

\begin{document}

\vspace*{1.4cm}

\normalsize \centerline{\Large \bf ON THE WELL-POSEDNESS
OF A CERTAIN MODEL}

\medskip

\centerline{\Large \bf WITH MIXED DIFFUSION APPEARING IN}

\medskip

\centerline{\Large \bf   THE MATHEMATICAL BIOLOGY}

\vspace*{1cm}

\centerline{\bf Messoud Efendiev$^{1,2}$, Vitali Vougalter$^{3 \ *}$}

\bigskip

\centerline{$^1$ Helmholtz Zentrum M\"unchen, Institut f\"ur Computational
Biology, Ingolst\"adter Landstrasse 1}

\centerline{Neuherberg, 85764, Germany}

\centerline{e-mail: messoud.efendiyev@helmholtz-muenchen.de}

\centerline{$^2$ Azerbaijan University of Architecture and Construction, Baku, Azerbaijan}

\centerline{e-mail: messoud.efendiyev@gmail.com}

\bigskip

\centerline{$^{3 \ *}$  Department of Mathematics, University
of Toronto}

\centerline{Toronto, Ontario, M5S 2E4, Canada}

\centerline{ e-mail: vitali@math.toronto.edu}

\medskip


\vspace*{0.25cm}

\noindent {\bf Abstract:}
The work deals with establishing the global well-posedness in the space
$W^{(1,2),2}({\mathbb R}^{2}\times {\mathbb R}^{+})$ of the
integro-differential equation  in the space of two dimensions in the situation when the
 diffusion term involves the standard Laplace operator in one direction and the fractional 
Laplacian in another direction in
the presence of the drift term. Such model is relevant to the cell
population dynamics in the Mathematical Biology.
The proof is based on a fixed point technique.

\vspace*{0.25cm}



\vspace*{0.5cm}

\bigskip

\bigskip


\setcounter{section}{1}

\centerline{\bf 1. Introduction}

\medskip

\noindent
The present article is devoted to the global well-posedness
of the nonlocal reaction-diffusion equation with the constants
$\displaystyle{0<s<1, \ a\geq 0}$, the constant vector $b:=(b_{1}, b_{2})\in {\mathbb R}^{2}$ of the length
$|b|_{{\mathbb R}^{2}}:=\sqrt{b_{1}^{2}+b_{2}^{2}}$, namely
\begin{equation}
\label{h}
\frac{\partial u}{\partial t} =\frac{{\partial^{2}}u}{\partial {x_{1}}^{2}}
-\Bigg(-\frac{\partial^{2}}{\partial x_{2}^{2}}\Bigg)^{s} u +
b.\nabla u+au+
\int_{{\mathbb R}^{2}}G(x-y)F(u(y,t), y)dy
\end{equation}
relevant to the cell population dynamics.  Here $x=(x_{1}, x_{2})\in {\mathbb R}^{2}, \ y=(y_{1}, y_{2})\in {\mathbb R}^{2}$.
The "dot" symbol in (\ref{h}) and below stands for the scalar product of two 
vectors in ${\mathbb R}^{2}$.
We assume here that the initial
condition for (\ref{h}) is
\begin{equation}
\label{ic}
u(x,0)=u_{0}(x)\in H^{2}({\mathbb R}^{2}).
\end{equation}  
Note that the existence of stationary solutions of the problem analogous to
(\ref{h}) without the transport term and for $a=0$ was established in ~\cite{EV210}.
The necessary conditions for preserving the nonnegative cone for a system of parabolic equations in the case of the mixed diffusion were derived in ~\cite{EV21}.
The stationary situation in the square with periodic boundary conditions and 
$\displaystyle{s=\frac{1}{2}}$ was covered in ~\cite{EV24}. The argument in ~\cite{EV24}
was performed in a constrained Sobolev space, which is a Hilbert space as well (see e.g.
Chapter 2.1 of ~\cite{HS96}). The article ~\cite{GHLP22} deals with the 
local and global existence for nonlocal multispecies advection-diffusion models.
The global well-posedness of the integro-differential equation analogous to (\ref{h}) considered
on the real line with  the fractional Laplacian in the diffusion term in the context of the anomlous diffusion was
established in ~\cite{EV25}. The anomalous diffusion can be described as a random process of particle motion 
characterized by the probability density distribution of jump length. The moments of this density distribution 
are finite in the case of the normal diffusion, but this is  not the case for the anomalous diffusion. 
The asymptotic behavior at infinity of the probability density function determines the value of the power of the Laplacian
(see ~\cite{MK00}).
The situation of the bi-Laplacian in the diffusion term was discussed in
~\cite{EV251} and a model with two kernels was covered in ~\cite{EV26}.
The stationary case of the normal diffusion and transport was considered in
~\cite{EV20} and the problem with drift and superdiffusion was solved in ~\cite{EV240}.
The work ~\cite{VV21} deals with the existence of
stationary solutions of the integro-differential equation on the real line containing the
fractional Laplacian raised to the power
$\displaystyle{0<s<\frac{1}{4}}$, the transport and the cell influx/efflux.
Existence of solutions for certain nonlinear equations with mixed local and nonlocal operators
was studied in ~\cite{I26}. Stability results for evolution equations involving superposition operators of mixed fractional
orders were derived in ~\cite{KK25}. The regularity and maximum principles for mixed local and nonlocal elliptic operators
were covered in ~\cite{BDVV22}. Certain nonlinear problems for the superposition of fractional operators with Neumann boundary
conditions were solved in ~\cite{DPSV25}.
The article ~\cite{ABVV10} is devoted to the 
spatial structures and generalized travelling waves for an
integro-differential equation. 
Emergence and propagation of patterns in nonlocal reaction-
diffusion equations arising in the theory of speciation and involving the
drift term were discussed in ~\cite{VV130}.  Pattern and waves for a model in population
dynamics with nonlocal consumption of resources were covered in ~\cite{GVA06}.
In ~\cite{BNPR09} the authors established the existence of steady states and travelling waves for the
non-local Fisher-KPP equation. The speed of propagation for KPP type
problems in the periodic framework was estimated in ~\cite{BHN05}. Significant applications to the theory of
reaction-diffusion problems with non-Fredholm operators were developed in
~\cite{DMV05}, ~\cite{DMV08}. The work ~\cite{GK18} deals with the entropy method for generalized Poisson-Nernst-
Planck equations. In ~\cite{DF06} the authors  demonstrated the exponential decay toward equilibrium via entropy methods for reaction-diffusion problems.  Lower estimate of the attractor dimension for a chemotaxis growth system was obtained in
~\cite{ATEYM06}. Quasilinear elliptic equations on half- and quarter-spaces were discussed in ~\cite{DDE13}.
Fredholm structures, topological invariants and applications were covered in ~\cite{E09}. The work ~\cite{E10} 
deals with the development of the theory of finite and infinite dimensional attractors for evolution equations of mathematical physics.
In ~\cite{E13} the author studied  evolution equations arising in the modelling of life sciences.
 Attractors for degenerate parabolic type equations were treated in ~\cite{E131}. The article ~\cite{EOE20} is devoted to the
mathematical analysis of a PDE–ODE coupled model of mitochondrial swelling with degenerate calcium ion diffusion.
In ~\cite{EP07} the authors investigated the large time behavior of solutions of fourth order parabolic equations and
$\epsilon$-entropy of their attractors.

The space variable $x$ in the present work is correspondent to the cell genotype,
$u(x,t)$ denotes the cell density as a function of the genotype and time.
The right side of (\ref{h}) describes the evolution of the cell density by means of
the cell proliferation, mutations and transport.
The mixed diffusion term here corresponds to the change of genotype
due to the small random mutations, and the integral nonlinear production term describes large mutations.
The function $F(u, x)$ stands for the rate of the cell birth, which depends on $u$ and $x$
(density dependent proliferation), and the kernel $G(x-y)$ designates
the proportion of newly born cells changing their genotype from $y$ to $x$.
Let us assume that it depends on the distance between the genotypes.

The operator involved in the mixed diffusion term of equation (\ref{h}) 
\begin{equation}
\label{ls}
L_{s}:=-\frac{{\partial^{2}}}{\partial {x_{1}}^{2}}
+\Bigg(-\frac{\partial^{2}}{\partial x_{2}^{2}}\Bigg)^{s}: H^{2}({\mathbb R}^{2})\to  L^{2}({\mathbb R}^{2}), \quad 0<s<1
\end{equation}
is defined via the spectral calculus. It is the pseudo-differential
operator with the symbol $p_{1}^{2}+|p_{2}|^{2s}$, such that
$$
L_{s}\phi(x)=\frac{1}{2\pi}
\int_{{\mathbb R}^{2}}(p_{1}^{2}+|p_{2}|^{2s}){\widehat \phi}(p)e^{ip.x}dp, \quad
\phi(x)\in H^{2}({\mathbb R}^{2})
$$
with the Sobolev space introduced in (\ref{ss}).
The standard Fourier transform in this context is given by
\begin{equation}
\label{ft}
\widehat{\phi}(p)=\frac{1}{2\pi}\int_{{\mathbb R}^{2}}\phi(x)e^{-ip.x}dx,
\quad p\in {\mathbb R}^{2}.
\end{equation}
Clearly, the estimate from above
\begin{equation}
\label{fub}  
\|\widehat{\phi}(p)\|_{L^{\infty}({\mathbb R}^{2})}\leq \frac{1}{2\pi}
\|\phi(x)\|_{L^{1}({\mathbb R}^{2})}
\end{equation}
is valid (see e.g. ~\cite{LL97}). Evidently, (\ref{fub}) gives us
\begin{equation}
\label{fub1}  
\|p^{2}\widehat{\phi}(p)\|_{L^{\infty}({\mathbb R}^{2})}\leq \frac{1}{2\pi}\|\Delta_{x} \phi\|_{L^{1}({\mathbb R}^{2})}.
\end{equation}
Here and below $\Delta_{x}$ stands for the Laplace operator with respect to the $x$ variable.
 The novelty of the present article is that in the diffusion term 
we add the standard minus Laplace operator in the $x_{1}$ variable to the negative  Laplacian in  $x_{2}$
raised to a fractional power.  This model is new and not much is understood about it, especially in 
the context of the integro-differential problems. The difficulty we have to overcome is that this 
model becomes anisotropic and it is more technical to derive the desired estimates when working  with it. 
Let us suppose that the conditions below on the integral kernel contained in 
equation (\ref{h}) are satisfied.

\medskip

\noindent
{\bf Assumption 1.1.}  {\it Let $G(x): {\mathbb R}^{2}\to {\mathbb R}$ be
nontrivial, such that}
$$
G(x), \ \Delta_{x} G(x)\in L^{1}({\mathbb R}^{2}).
$$     

\medskip

This allows us to define the auxiliary quantity
\begin{equation}
\label{q}  
\sigma:=\sqrt{\|G(x)\|_{L^{1}({\mathbb R}^{2})}^{2}+
\|\Delta_{x} G(x)\|_{L^{1}({\mathbb R}^{2})}^{2}}.
\end{equation}
Obviously, $0<\sigma<\infty$.
  
\medskip

From the point of view of the applications, the space dimension is
not limited to $d=2$ since the space variable corresponds to the cell
genotype but not to the usual physical space.
We are using the Sobolev space
\begin{equation}
\label{ss}  
H^{2}({\mathbb R}^{2}):=
\Big\{\phi(x):{\mathbb R}^{2}\to {\mathbb {\mathbb R}} \ | \
\phi(x)\in L^{2}({\mathbb R}^{2}), \ \Delta_{x} \phi(x) \in
L^{2}({\mathbb R}^{2}) \Big \}.
\end{equation}
It is equipped with the norm
\begin{equation}
\label{n}
\|\phi\|_{H^{2}({\mathbb R}^{2})}^{2}:=\|\phi\|_{L^{2}({\mathbb R}^{2})}^{2}+
\| \Delta_{x} \phi\|_{L^{2}({\mathbb R}^{2})}^{2}.
\end{equation}
Let us establish the global well-posedness for problem (\ref{h}), (\ref{ic}).
For that purpose, we will use the function space
$$
W^{(1, 2), 2}({\mathbb R}^{2}\times [0, T]):=
$$
\begin{equation}
\label{122}
\Big\{u(x,t): {\mathbb R}^{2}\times [0, T]\to {\mathbb R} \ \Big|
\ u(x,t), \ \Delta_{x}u, \
\frac{\partial u}{\partial t}
\in L^{2}({\mathbb R}^{2}\times [0, T]) \Big\}.
\end{equation}  
Hence,
$$
\|u(x,t)\|_{W^{(1, 2), 2}({\mathbb R}^{2}\times [0, T])}^{2}:=
$$
\begin{equation}
\label{122n}
\Big\|\frac{\partial u}{\partial t}\Big\|_{L^{2}({\mathbb R}^{2}\times [0, T])}^{2}+
\|\Delta_{x}u\|_{L^{2}({\mathbb R}^{2}\times [0, T])}^{2}+
\|u\|_{L^{2}({\mathbb R}^{2}\times [0, T])}^{2},
\end{equation}  
where $T>0$ and
$$
\|u\|_{L^{2}({\mathbb R}^{2}\times [0, T])}^{2}:=\int_{0}^{T}\int_{{\mathbb R}^{2}}
|u(x,t)|^{2}dxdt.
$$
In the present work we will also use another norm
$$
\|u(x,t)\|_{L^{2}({\mathbb R}^{2})}^{2}:=\int_{{\mathbb R}^{2}}|u(x,t)|^{2}dx.
$$

\medskip

\noindent
{\bf Assumption 1.2.} {\it Function
$F(u, x): {\mathbb R}\times{\mathbb R}^{2}\to {\mathbb R}$ is satisfying the
Caratheodory condition (see ~\cite{K64}), so that
\begin{equation}
\label{lub}
|F(u, x)|\leq k|u|+h(x) \quad for \quad u\in {\mathbb R}, \quad x\in {\mathbb R}^{2}
\end{equation}  
with a constant $k>0$ and
$h(x): {\mathbb R}^{2}\to {\mathbb R}^{+}, \ h(x)\in L^{2}({\mathbb R}^{2})$.
Moreover, it is a Lipschitz continuous function, such that
\begin{equation}
\label{lip}
|F(u_{1}, x)-F(u_{2}, x)|\leq l|u_{1}-u_{2}| \quad for \quad any \quad
u_{1, 2}\in {\mathbb R}, \quad x\in {\mathbb R}^{2}
\end{equation}
with a constant $l>0$.}

\medskip

In the article ${\mathbb R}^{+}$ designates the nonnegative semi-axis.
The solvability of a local elliptic problem in a bounded domain in
${\mathbb R}^{N}$ was covered in ~\cite{BO86}. The nonlinear function involved there
was allowed to have a sublinear growth.
We apply the standard Fourier transform (\ref{ft}) to both sides of
(\ref{h}), (\ref{ic}) and arrive at
\begin{equation}
\label{hf}
\frac{\partial \widehat{u}}{\partial t}=[-p_{1}^{2}-|p_{2}|^{2s}+ib.p+a]\widehat{u}+
2\pi \widehat{G}(p)\widehat{f}_{u}(p,t),
\end{equation}
\begin{equation}
\label{icf}
\widehat{u(x,0)}(p)=\widehat{u_{0}}(p).  
\end{equation}
Here and below $\widehat{f}_{u}(p,t)$ will denote the Fourier image
of $F(u(x,t), x)$. Evidently, we have
$$
u(x,t)=\frac{1}{2\pi}\int_{{\mathbb R}^{2}}\widehat{u}(p,t)e^{ip.x}dp,
\quad 
\frac{\partial u}{\partial t}=\frac{1}{2\pi}\int_{{\mathbb R}^{2}}
\frac{\partial \widehat{u}(p,t)}{\partial t}e^{ip.x}dp,
$$
where $x\in {\mathbb R}^{2}, \ t\geq 0$. According to the Duhamel's principle,
we can reformulate (\ref{hf}), (\ref{icf}) as
$$
\widehat{u}(p,t)=
$$
\begin{equation}
\label{duh}
e^{t\{-p_{1}^{2}-|p_{2}|^{2s}+ib.p+a\}}\widehat{u_{0}}(p)+\int_{0}^{t}
e^{(t-{\tilde s})\{-p_{1}^{2}-|p_{2}|^{2s}+ib.p+a\}}2\pi\widehat{G}(p)\widehat{f}_{u}(p, {\tilde s})d{\tilde s}.
\end{equation}
Let us consider the auxiliary equation related to problem (\ref{duh}), namely
$$
\widehat{u}(p,t)=
$$
\begin{equation}
\label{aux}
e^{t\{-p_{1}^{2}-|p_{2}|^{2s}+ib.p+a\}}\widehat{u_{0}}(p)+\int_{0}^{t}
e^{(t-{\tilde s})\{-p_{1}^{2}-|p_{2}|^{2s}+ib.p+a\}}2\pi \widehat{G}(p)\widehat{f}_{v}(p, {\tilde s})d{\tilde s}.
\end{equation}
Here $\widehat{f}_{v}(p, {\tilde s})$ stands for the Fourier image of
$F(v(x, {\tilde s}), x)$ under transform (\ref{ft}) and $v(x,t)\in  W^{(1, 2), 2}({\mathbb R}^{2}\times [0, T])$
is arbitrary.

We introduce the operator $T_{a, b}$, such that $u =T_{a, b}v$, with $u$ 
satisfying (\ref{aux}). The main proposition of the work is as follows.

\bigskip

\noindent
{\bf Theorem 1.3.} {\it Let Assumptions 1.1 and 1.2 be valid and
\begin{equation}
\label{qlt}
\sigma l\sqrt{T^{2}e^{2aT}+2\alpha_{a, b, T}}<1
\end{equation}
with the constant $\sigma$ defined in (\ref{q}), the Lipschitz constant
$l$ introduced in (\ref{lip})  and $\alpha_{a, b, T}$ given by (\ref{alf}).
Then problem (\ref{aux}) defines the map
$T_{a, b}: W^{(1, 2), 2}({\mathbb R}^{2}\times [0, T])\to W^{(1, 2), 2}
({\mathbb R}^{2}\times [0, T])$, which is a strict contraction.
The unique fixed point $w(x,t)$ of this map $T_{a, b}$ is the only solution of
(\ref{h}), (\ref{ic}) in
$W^{(1, 2), 2}({\mathbb R}^{2}\times [0, T])$.}

\medskip

The final result of our article is devoted to the global well-posedness for
our equation. 

\medskip

\noindent
{\bf Corollary 1.4.} {\it Let the assumptions of Theorem 1.3 hold.
Then problem (\ref{h}), (\ref{ic}) possesses a unique solution   
$w(x,t)\in W^{(1, 2), 2}({\mathbb R}^{2}\times {\mathbb R}^{+})$. Such solution 
does not vanish identically for $x\in {\mathbb R}^{2}$ and $t\in {\mathbb R}^{+}$ 
provided the intersection of supports of the Fourier images of the functions
$\hbox{supp}\widehat{F(0, x)}\cap \hbox{supp}\widehat{G}$ is a set of
nonzero Lebesgue measure in ${\mathbb R}^{2}$.}

\medskip

Let us proceed to the proof of our main statement.

\bigskip


\setcounter{section}{2}
\setcounter{equation}{0}

\centerline{\bf 2. The well-posedness of the model}

\bigskip

\noindent
{\it Proof of Theorem 1.3.} We choose arbitrarily
$v(x,t)\in W^{(1, 2), 2}({\mathbb R}^{2}\times [0, T])$.
Let us demonstrate that the first term in the right side of (\ref{aux}) is
contained in $L^{2}({\mathbb R}^{2}\times [0, T])$. Clearly,
$$
\|e^{t\{-p_{1}^{2}-|p_{2}|^{2s}+ib.p+a\}}\widehat{u_{0}}(p)\|_{L^{2}({\mathbb R}^{2})}^{2}=
\int_{{\mathbb R}^{2}}e^{-2tp_{1}^{2}}e^{-2t|p_{2}|^{2s}}e^{2at}|\widehat{u_{0}}(p)|^{2}dp\leq
$$
$$
e^{2at}\|u_{0}\|_{L^{2}({\mathbb R}^{2})}^{2}.
$$
Then
$$
\|e^{t\{-p_{1}^{2}-|p_{2}|^{2s}+ib.p+a\}}\widehat{u_{0}}(p)\|_{L^{2}({\mathbb R}^{2}\times [0, T])}^{2}=
\int_{0}^{T}\|e^{t\{-p_{1}^{2}-|p_{2}|^{2s}+ib.p+a\}}\widehat{u_{0}}(p)\|_{L^{2}({\mathbb R}^{2})}^{2}dt
\leq
$$
$$
\int_{0}^{T}e^{2at}\|u_{0}\|_{L^{2}({\mathbb R}^{2})}^{2}dt.
$$
Obviously, this is equal to
$\displaystyle{\frac{e^{2aT}-1}{2a}\|u_{0}\|_{L^{2}({\mathbb R}^{2})}^{2}}$
if $a>0$ and
$T\|u_{0}\|_{L^{2}({\mathbb R}^{2})}^{2}$ for $a=0$. This means that
\begin{equation}
\label{u0l2}  
e^{t\{-p_{1}^{2}-|p_{2}|^{2s}+ib.p+a\}}\widehat{u_{0}}(p)\in L^{2}({\mathbb R}^{2}\times [0, T]).
\end{equation}
We obtain the upper bound on the norm of the second term in the right side of (\ref{aux}) as
$$
\Big\|\int_{0}^{t}
e^{(t-{\tilde s})\{-p_{1}^{2}-|p_{2}|^{2s}+ib.p+a\}}2\pi \widehat{G}(p)\widehat{f}_{v}(p, {\tilde s})d{\tilde s}\Big\|_
{L^{2}({\mathbb R}^{2})}\leq
$$
$$  
\int_{0}^{t}\Big\|
e^{(t-{\tilde s})\{-p_{1}^{2}-|p_{2}|^{2s}+ib.p+a\}}2\pi \widehat{G}(p)\widehat{f}_{v}(p,{\tilde s})\Big\|_
{L^{2}({\mathbb R}^{2})}d{\tilde s}.
$$
Evidently,
$$  
\Big\|
e^{(t-{\tilde s})\{-p_{1}^{2}-|p_{2}|^{2s}+ib.p+a\}}2\pi\widehat{G}(p)\widehat{f}_{v}(p,{\tilde s})\Big\|_
{L^{2}({\mathbb R}^{2})}^{2}=
$$
\begin{equation}
\label{intgfv}  
\int_{{\mathbb R}^{2}}e^{-2(t-{\tilde s})p_{1}^{2}}e^{-2(t-{\tilde s})|p_{2}|^{2s}}e^{2a(t-{\tilde s})}(2\pi)^{2}|\widehat{G}(p)|^{2}
|\widehat{f}_{v}(p, {\tilde s})|^{2}dp.
\end{equation}
Let us apply inequality (\ref{fub}) to derive the estimate from above on the right side of
(\ref{intgfv}) as
$$
e^{2a(t-{\tilde s})}(2\pi)^{2}\|\widehat{G}(p)\|_{L^{\infty}({\mathbb R}^{2})}^{2}
\|F(v(x, {\tilde s}),x)\|_{L^{2}({\mathbb R}^{2})}^{2}\leq
$$
$$
e^{2aT}\|G(x)\|_{L^{1}({\mathbb R}^{2})}^{2}
\|F(v(x, {\tilde s}),x)\|_{L^{2}({\mathbb R}^{2})}^{2}.
$$
Then
$$  
\Big\|
e^{(t-{\tilde s})\{-p_{1}^{2}-|p_{2}|^{2s}+ib.p+a\}}2\pi \widehat{G}(p)\widehat{f}_{v}(p, {\tilde s})\Big\|_
{L^{2}({\mathbb R}^{2})}\leq
$$
$$
e^{aT}\|G(x)\|_{L^{1}({\mathbb R}^{2})}
\|F(v(x, {\tilde s}),x)\|_{L^{2}({\mathbb R}^{2})}.
$$
According to (\ref{lub}), we have
\begin{equation}
\label{lubs}  
\|F(v(x, {\tilde s)},x)\|_{L^{2}({\mathbb R}^{2})}\leq k\|v(x, {\tilde s})\|_{L^{2}({\mathbb R}^{2})}+\|h(x)\|_
{L^{2}({\mathbb R}^{2})}.
\end{equation}
This yields
$$  
\Big\|
e^{(t-{\tilde s})\{-p_{1}^{2}-|p_{2}|^{2s}+ib.p+a\}}2\pi \widehat{G}(p)\widehat{f}_{v}(p, {\tilde s})\Big\|_
{L^{2}({\mathbb R}^{2})}\leq
$$
$$
e^{aT}\|G(x)\|_{L^{1}({\mathbb R}^{2})}
\{k\|v(x, {\tilde s})\|_{L^{2}({\mathbb R}^{2})}+\|h(x)\|_{L^{2}({\mathbb R}^{2})}\},
$$
so that
$$  
\Big\|\int_{0}^{t}
e^{(t-{\tilde s})\{-p_{1}^{2}-|p_{2}|^{2s}+ib.p+a\}}2\pi \widehat{G}(p)\widehat{f}_{v}(p, {\tilde s})d {\tilde s}\Big\|_
{L^{2}({\mathbb R}^{2})}\leq
$$
$$
ke^{aT}\|G(x)\|_{L^{1}({\mathbb R}^{2})}\int_{0}^{T}\|v(x, {\tilde s})\|_{L^{2}({\mathbb R}^{2})}d {\tilde s}+
Te^{aT}\|G(x)\|_{L^{1}({\mathbb R}^{2})}\|h(x)\|_{L^{2}({\mathbb R}^{2})}.
$$
Let us recall the Schwarz inequality. Hence,
\begin{equation}
\label{sch}  
\int_{0}^{T}\|v(x,{\tilde s})\|_{L^{2}({\mathbb R}^{2})}d {\tilde s}\leq
\sqrt{\int_{0}^{T}\|v(x, {\tilde s})\|_{L^{2}({\mathbb R}^{2})}^{2}d {\tilde s}}\sqrt{T}.
\end{equation}
This gives us
$$
\Big\|\int_{0}^{t}
e^{(t- {\tilde s})\{-p_{1}^{2}-|p_{2}|^{2s}+ib.p+a\}}2\pi \widehat{G}(p)\widehat{f}_{v}(p, {\tilde s})d {\tilde s}\Big\|_
{L^{2}({\mathbb R}^{2})}^{2}\leq
$$
$$
e^{2aT}\|G(x)\|_{L^{1}({\mathbb R}^{2})}^{2}
\{k\sqrt{T}\|v(x, {\tilde s})\|_{L^{2}({\mathbb R}^{2}\times [0,T])}+T\|h(x)\|_{L^{2}({\mathbb R}^{2})}\}^{2}.
$$
We derive the upper bound on the norm as
$$
\Big\|\int_{0}^{t}
e^{(t-{\tilde s})\{-p_{1}^{2}-|p_{2}|^{2s}+ib.p+a\}}2\pi  \widehat{G}(p)\widehat{f}_{v}(p, {\tilde s})d {\tilde s}\Big\|_
{L^{2}({\mathbb R}^{2}\times [0,T])}^{2}=
$$
$$
\int_{0}^{T}\Big\|\int_{0}^{t}
e^{(t-{\tilde s})\{-p_{1}^{2}-|p_{2}|^{2s}+ib.p+a\}}2\pi \widehat{G}(p)\widehat{f}_{v}(p,{\tilde s})d {\tilde s}\Big\|_
{L^{2}({\mathbb R}^{2})}^{2}dt\leq 
$$
\begin{equation}
\label{e2atgkvrtht}
e^{2aT}\|G(x)\|_{L^{1}({\mathbb R}^{2})}^{2}
\{k\|v(x, {\tilde s})\|_{L^{2}({\mathbb R}^{2}\times [0,T])}+\sqrt{T}\|h(x)\|_{L^{2}({\mathbb R}^{2})}\}^{2}
T^{2}.
\end{equation}
Note that the right side of (\ref{e2atgkvrtht}) is finite 
under the stated assumptions for
$v(x,t)\in W^{(1, 2), 2}({\mathbb R}^{2}\times [0, T])$. This proves that
\begin{equation}
\label{0tetsgvl2}  
\int_{0}^{t}
e^{(t-{\tilde s})\{-p_{1}^{2}-|p_{2}|^{2s}+ib.p+a\}}2\pi \widehat{G}(p)\widehat{f}_{v}(p, {\tilde s})d {\tilde s}\in
L^{2}({\mathbb R}^{2}\times [0,T]).
\end{equation}
By virtue of (\ref{aux}), (\ref{u0l2}) and (\ref{0tetsgvl2}), we obtain
\begin{equation}
\label{upt12}  
\widehat{u}(p,t)\in L^{2}({\mathbb R}^{2}\times [0,T]).
\end{equation}
Thus,
\begin{equation}
\label{uxtl2}  
u(x,t)\in L^{2}({\mathbb R}^{2}\times [0,T]).
\end{equation}
Let us recall (\ref{aux}). Hence,
$$
p^{2}\widehat{u}(p,t)=
$$
\begin{equation}
\label{aux2}
e^{t\{-p_{1}^{2}-|p_{2}|^{2s}+ib.p+a\}}p^{2}\widehat{u_{0}}(p)+\int_{0}^{t}
e^{(t-{\tilde s})\{-p_{1}^{2}-|p_{2}|^{2s}+ib.p+a\}}2\pi p^{2}\widehat{G}(p)\widehat{f}_{v}(p, {\tilde s})d {\tilde s}.
\end{equation}
We estimate the norm of the first term in the right side of (\ref{aux2}) as
$$
\|e^{t\{-p_{1}^{2}-|p_{2}|^{2s}+ib.p+a\}}p^{2}\widehat{u_{0}}(p)\|_{L^{2}({\mathbb R}^{2}\times [0,T])}^{2}=
$$
$$
\int_{0}^{T}\int_{{\mathbb R}^{2}}e^{-2tp_{1}^{2}}e^{-2t|p_{2}|^{2s}}e^{2at}|p^{2}\widehat{u_{0}}(p)|^{2}dp
dt\leq
\int_{0}^{T}\int_{{\mathbb R}^{2}}e^{2at}|p^{2}\widehat{u_{0}}(p)|^{2}dpdt.
$$
Clearly, this equals to
$\displaystyle
{\frac{e^{2aT}-1}{2a}\|\Delta_{x} u_{0}\|_{L^{2}({\mathbb R}^{2})}^{2}}$ for
$a>0$ and
$\displaystyle{T\|\Delta_{x} u_{0}\|_{L^{2}({\mathbb R}^{2})}^{2}}$
for $a=0$. Therefore,
\begin{equation}
\label{p2u0hpl2}
e^{t\{-p_{1}^{2}-|p_{2}|^{2s}+ib.p+a\}}p^{2}\widehat{u_{0}}(p)\in L^{2}({\mathbb R}^{2}\times [0,T]).
\end{equation}  
We proceed to analyzing the second term in the right side of
(\ref{aux2}). Obviously,
$$
\Big\|\int_{0}^{t}
e^{(t-{\tilde s})\{-p_{1}^{2}-|p_{2}|^{2s}+ib.p+a\}}2\pi p^{2}\widehat{G}(p)\widehat{f}_{v}(p, {\tilde s})d {\tilde s}
\Big\|_{L^{2}({\mathbb R}^{2})}\leq 
$$
$$
\int_{0}^{t}\Big\|
e^{(t-{\tilde s})\{-p_{1}^{2}-|p_{2}|^{2s}+ib.p+a\}}2\pi p^{2}\widehat{G}(p)\widehat{f}_{v}(p, {\tilde s})
\Big\|_{L^{2}({\mathbb R}^{2})}d {\tilde s}.
$$
We easily derive
$$
\Big\|e^{(t-{\tilde s})\{-p_{1}^{2}-|p_{2}|^{2s}+ib.p+a\}}2\pi p^{2}\widehat{G}(p)\widehat{f}_{v}(p, {\tilde s})
\Big\|_{L^{2}({\mathbb R}^{2})}^{2}=
$$
\begin{equation}
\label{etsp2gpl2}
\int_{{\mathbb R}^{2}}
e^{-2(t-{\tilde s})p_{1}^{2}}e^{-2(t- {\tilde s})|p_{2}|^{2s}}e^{2a(t-{\tilde s})}(2\pi)^{2}|p^{2}\widehat{G}(p)|^{2}|\widehat{f}_{v}(p, {\tilde s})|^{2}
dp.
\end{equation}  
The right side of (\ref{etsp2gpl2}) can be bounded above using (\ref{fub1})
as
$$
(2\pi)^{2} e^{2aT}\|p^{2}\widehat{G}(p)\|_{L^{\infty}({\mathbb R}^{2})}^{2}
\int_{{\mathbb R}^{2}}|\widehat{f}_{v}(p, {\tilde s})|^{2}dp\leq
$$
$$
e^{2aT}\|\Delta_{x} G(x) \|_{L^{1}({\mathbb R}^{2})}^{2}
\|F(v(x, {\tilde s}),x)\|_{L^{2}({\mathbb R}^{2})}^{2}.
$$
Let us use estimate (\ref{lubs}) to obtain the inequality
$$
\|e^{(t-{\tilde s})\{-p_{1}^{2}-|p_{2}|^{2s}+ib.p+a\}}2\pi p^{2}\widehat{G}(p)\widehat{f}_{v}(p, {\tilde s})
\Big\|_{L^{2}({\mathbb R}^{2})}\leq
$$
$$
e^{aT}\|\Delta_{x} G(x) \|_{L^{1}({\mathbb R}^{2})}\{k\|v(x, {\tilde s})\|_
{L^{2}({\mathbb R}^{2})}+\|h(x)\|_{L^{2}({\mathbb R}^{2})}\}.
$$
Then
$$
\Big\|\int_{0}^{t}e^{(t-{\tilde s})\{-p_{1}^{2}-|p_{2}|^{2s}+ib.p+a\}}2\pi p^{2}\widehat{G}(p)
\widehat{f}_{v}(p, {\tilde s})d {\tilde s}\Big\|_{L^{2}({\mathbb R}^{2})}\leq
$$
$$
ke^{aT}\|\Delta_{x} G(x)\|_{L^{1}({\mathbb R}^{2})}\int_{0}^{T}\|v(x, {\tilde s})\|_
{L^{2}({\mathbb R}^{2})}d {\tilde s}+
Te^{aT}\|\Delta_{x} G(x)\|_{L^{1}({\mathbb R}^{2})}\|h(x)\|_{L^{2}({\mathbb R}^{2})}.
$$
Recall upper bound (\ref{sch}). We arrive at
$$
\Big\|\int_{0}^{t}e^{(t-{\tilde s})\{-p_{1}^{2}-|p_{2}|^{2s}+ib.p+a\}}2\pi p^{2}\widehat{G}(p)
\widehat{f}_{v}(p, {\tilde s})d {\tilde s}\Big\|_{L^{2}({\mathbb R}^{2})}^{2}\leq
$$
$$
e^{2aT}\|\Delta_{x} G(x)\|_{L^{1}({\mathbb R}^{2})}^{2}
\{k\|v(x, {\tilde s})\|_{L^{2}({\mathbb R}^{2}\times [0,T])}\sqrt{T}+
\|h(x)\|_{L^{2}({\mathbb R}^{2})}T\}^{2}.
$$
This means that
$$
\Big\|\int_{0}^{t}e^{(t-{\tilde s})\{-p_{1}^{2}-|p_{2}|^{2s}+ib.p+a\}}2\pi p^{2}\widehat{G}(p)
\widehat{f}_{v}(p, {\tilde s})d {\tilde s}\Big\|_{L^{2}({\mathbb R}^{2}\times [0, T])}^{2}\leq
$$
\begin{equation}
\label{e2atdgkvht2}
e^{2aT}\|\Delta_{x}G(x)\|_{L^{1}({\mathbb R}^{2})}^{2}
\{k\|v(x, {\tilde s})\|_{L^{2}({\mathbb R}^{2}\times [0,T])}+
\|h(x)\|_{L^{2}({\mathbb R}^{2})}\sqrt{T}\}^{2}T^{2}.
\end{equation}
We use the given conditions and the fact that
$v(x,t)\in W^{(1, 2), 2}({\mathbb R}^{2}\times [0, T])$. Thus, the right side of
(\ref{e2atdgkvht2}) is finite. This proves that
\begin{equation}
\label{int0tp2Gpfv}  
\int_{0}^{t}e^{(t-{\tilde s})\{-p_{1}^{2}-|p_{2}|^{2s}+ib.p+a\}}2\pi p^{2}\widehat{G}(p)
\widehat{f}_{v}(p, {\tilde s})d {\tilde s}\in L^{2}({\mathbb R}^{2}\times [0, T]).
\end{equation}
According to  (\ref{aux2}), (\ref{p2u0hpl2}) and (\ref{int0tp2Gpfv}), we have
\begin{equation}
\label{up2t12}  
p^{2}\widehat{u}(p,t)\in L^{2}({\mathbb R}^{2}\times [0, T]).
\end{equation}  
Therefore,
\begin{equation}
\label{d2udx2l2}
\Delta_{x} u \in L^{2}({\mathbb R}^{2}\times [0, T]).
\end{equation}  
It follows easily from (\ref{aux}) that
\begin{equation}
\label{duhdt}
\frac{\partial \widehat{u}}{\partial t}=\{-p_{1}^{2}-|p_{2}|^{2s}+ib.p+a\}\widehat{u}(p,t)+
2\pi \widehat{G}(p)\widehat{f}_{v}(p,t).
\end{equation}
Let us recall (\ref{upt12}). Hence,
\begin{equation}
\label{aupt12}  
a\widehat{u}(p,t)\in L^{2}({\mathbb R}^{2}\times [0,T]).
\end{equation}
We derive the estimate on the norm as
$$
\|ib.p\widehat{u}(p,t)\|_{L^{2}({\mathbb R}^{2}\times [0,T])}^{2}\leq
$$
$$ 
|b|_{{\mathbb R}^{2}}^{2}\int_{0}^{T}
\Big\{\int_{|p|\leq 1}p^{2}|\widehat{u}(p,t)|^{2}dp+
\int_{|p|>1}p^{2}|\widehat{u}(p,t)|^{2}dp\Big\}dt\leq
$$
\begin{equation}
\label{ibpuhn}
|b|_{{\mathbb R}^{2}}^{2}\{\|\widehat{u}(p,t)\|_{L^{2}({\mathbb R}^{2}\times [0,T])}^{2}+
\|p^{2}\widehat{u}(p,t)\|_{L^{2}({\mathbb R}^{2}\times [0,T])}^{2}\}.
\end{equation}
Clearly, the right side of (\ref{ibpuhn}) is finite due to
(\ref{upt12}) and (\ref{up2t12}). This means that
\begin{equation}
\label{bupt12} 
ib.p\widehat{u}(p,t)\in L^{2}({\mathbb R}^{2}\times [0,T]).
\end{equation}
Analogously,
$$
\|-|p_{2}|^{2s}\widehat{u}(p,t)\|_{L^{2}({\mathbb R}^{2}\times [0,T])}^{2}=
$$
$$
\int_{0}^{T}
\Big\{\int_{|p|\leq 1}|p_{2}|^{4s}|\widehat{u}(p,t)|^{2}dp+
\int_{|p|>1}|p_{2}|^{4s}|\widehat{u}(p,t)|^{2}dp\Big\}dt\leq
$$
$$
\|\widehat{u}(p,t)\|_{L^{2}({\mathbb R}^{2}\times [0,T])}^{2}+
\|p^{2}\widehat{u}(p,t)\|_{L^{2}({\mathbb R}^{2}\times [0,T])}^{2},
$$
which is finite via (\ref{upt12}) and (\ref{up2t12}). Thus
\begin{equation}
\label{uptalf12} 
-|p_{2}|^{2s}\widehat{u}(p,t)\in L^{2}({\mathbb R}^{2}\times [0,T]).
\end{equation}
Similarly,
$$
\|-p_{1}^{2}\widehat{u}(p,t)\|_{L^{2}({\mathbb R}^{2}\times [0,T])}^{2}=
$$
$$
\int_{0}^{T}
\Big\{\int_{|p|\leq 1}|p_{1}|^{4}|\widehat{u}(p,t)|^{2}dp+
\int_{|p|>1}|p_{1}|^{4}|\widehat{u}(p,t)|^{2}dp\Big\}dt\leq
$$
$$
\|\widehat{u}(p,t)\|_{L^{2}({\mathbb R}^{2}\times [0,T])}^{2}+
\|p^{2}\widehat{u}(p,t)\|_{L^{2}({\mathbb R}^{2}\times [0,T])}^{2}<\infty
$$
as above. This proves that
\begin{equation}
\label{mp12uhl2}
-p_{1}^{2}\widehat{u}(p,t)\in L^{2}({\mathbb R}^{2}\times [0,T]).
\end{equation}
Let us combine statements (\ref{aupt12}), (\ref{bupt12}), (\ref{uptalf12}) and (\ref{mp12uhl2}). We arrive at
\begin{equation}
\label{uptalfab12}
(-p_{1}^{2}-|p_{2}|^{2s}+ib.p+a)\widehat{u}(p,t)\in L^{2}({\mathbb R}^{2}\times [0,T]).
\end{equation}
We proceed to deriving the estimate from above on the norm of the
remaining term in the right side of (\ref{duhdt}) by means of (\ref{fub}) and 
(\ref{lubs}). Evidently,
$$
\|2\pi \widehat{G}(p)\widehat{f}_{v}(p,t)\|_
{L^{2}({\mathbb R}^{2}\times [0,T])}^{2}\leq (2\pi)^{2}
\|\widehat{G}(p)\|_{L^{\infty}({\mathbb R}^{2})}^{2}\int_{0}^{T}
\|F(v(x,t), x)\|_{L^{2}({\mathbb R}^{2})}^{2}dt\leq
$$
$$
\|G(x)\|_{L^{1}({\mathbb R}^{2})}^{2}\int_{0}^{T}(k\|v(x,t)\|_{L^{2}({\mathbb R}^{2})}+
\|h(x)\|_{L^{2}({\mathbb R}^{2})})^{2}dt\leq
$$
\begin{equation}
\label{2piGfvn}
\|G(x)\|_{L^{1}({\mathbb R}^{2})}^{2}\{2k^{2}\|v(x,t)\|_{L^{2}({\mathbb R}^{2}\times [0, T])}^{2}+
2\|h(x)\|_{L^{2}({\mathbb R}^{2})}^{2}T\}.
\end{equation}
Obviously, the right side of (\ref{2piGfvn}) is finite under our assumptions with 
$v(x,t)\in W^{(1, 2), 2}({\mathbb R}^{2}\times [0, T])$. Therefore,
\begin{equation}
\label{ghfvhpt}  
2\pi \widehat{G}(p)\widehat{f}_{v}(p,t)\in L^{2}({\mathbb R}^{2}\times [0,T]).
\end{equation}
Let us recall (\ref{duhdt}). By means of (\ref{uptalfab12}) along with
(\ref{ghfvhpt}), we obtain
$$
\frac{\partial \widehat{u}}{\partial t}\in L^{2}({\mathbb R}^{2}\times [0,T]).
$$
Hence,
\begin{equation}
\label{dudtl2}  
\frac{\partial u}{\partial t}\in L^{2}({\mathbb R}^{2}\times [0,T]).
\end{equation}
According to the definition of the norm (\ref{122n}) along with statements
(\ref{uxtl2}), (\ref{d2udx2l2}) and (\ref{dudtl2}), we derive that
for the function determined uniquely by equation (\ref{aux}), we have
$$
u(x,t)\in W^{(1, 2), 2}({\mathbb R}^{2}\times [0, T]).
$$
This means under the given conditions the auxiliary problem (\ref{aux}) defines the map
$$
T_{a, b}: W^{(1, 2), 2}({\mathbb R}^{2}\times [0, T])\to
W^{(1, 2), 2}({\mathbb R}^{2}\times [0, T]).
$$
Let us establish that under our assumptions such map is a
strict contraction. For this purpose, we choose arbitrarily
$v_{1, 2}(x,t)\in W^{(1, 2), 2}({\mathbb R}^{2}\times [0, T])$. According to the
reasoning above,
$u_{1, 2}:=T_{a, b}v_{1, 2}\in W^{(1, 2), 2}({\mathbb R}^{2}\times [0, T])$.
By virtue of (\ref{aux}),
$$
\widehat{u_{1}}(p,t)=
$$
\begin{equation}
\label{1aux}
e^{t\{-p_{1}^{2}-|p_{2}|^{2s}+ib.p+a\}}\widehat{u_{0}}(p)+\int_{0}^{t}
e^{(t-{\tilde s})\{-p_{1}^{2}-|p_{2}|^{2s}+ib.p+a\}}2\pi \widehat{G}(p)\widehat{f}_{v_{1}}(p, {\tilde s})d {\tilde s},
\end{equation}
$$
\widehat{u_{2}}(p,t)=
$$
\begin{equation}
\label{2aux}
e^{t\{-p_{1}^{2}-|p_{2}|^{2s}+ib.p+a\}}\widehat{u_{0}}(p)+\int_{0}^{t}
e^{(t-{\tilde s})\{-p_{1}^{2}-|p_{2}|^{2s}+ib.p+a\}}2\pi\widehat{G}(p)\widehat{f}_{v_{2}}(p, {\tilde s})d {\tilde s},
\end{equation}
where $\widehat{f}_{v_{j}}(p, {\tilde s})$ with $j=1, 2$ denotes the Fourier image of
$F(v_{j}(x, {\tilde s}), x)$ under transform (\ref{ft}). From the system of equations
(\ref{1aux}), (\ref{2aux}) we easily deduce that
$$
\widehat{u_{1}}(p,t)-\widehat{u_{2}}(p,t)=
$$
\begin{equation}
\label{u1u2hint0t}  
\int_{0}^{t}
e^{(t-{\tilde s})\{-p_{1}^{2}-|p_{2}|^{2s}+ib.p+a\}}2\pi \widehat{G}(p)
[\widehat{f}_{v_{1}}(p, {\tilde s})-\widehat{f}_{v_{2}}(p, {\tilde s})]d {\tilde s}.
\end{equation}
We derive the bound on the norm as
$$
\|\widehat{u_{1}}(p,t)-\widehat{u_{2}}(p,t)\|_{L^{2}({\mathbb R}^{2})}\leq
$$
\begin{equation}
\label{int0tetsf12}  
\int_{0}^{t}\|
e^{(t-{\tilde s})\{-p_{1}^{2}-|p_{2}|^{2s}+ib.p+a\}}2\pi\widehat{G}(p)
[\widehat{f}_{v_{1}}(p, {\tilde s)}-\widehat{f}_{v_{2}}(p, {\tilde s})]\|_{L^{2}({\mathbb R}^{2})}d {\tilde s}.
\end{equation}
By means of (\ref{fub}), we arrive at the estimate
$$
\|e^{(t-{\tilde s})\{-p_{1}^{2}-|p_{2}|^{2s}+ib.p+a\}}2\pi \widehat{G}(p)
[\widehat{f}_{v_{1}}(p, {\tilde s})-\widehat{f}_{v_{2}}(p, {\tilde s})]\|_{L^{2}({\mathbb R}^{2})}^{2}=
$$
$$
(2\pi)^{2} \int_{{\mathbb R}^{2}}e^{-2(t-{\tilde s})p_{1}^{2}}e^{-2(t-{\tilde s})|p_{2}|^{2s}}e^{2(t-{\tilde s})a}|\widehat{G}(p)|^{2}
|\widehat{f}_{v_{1}}(p,{\tilde s})-\widehat{f}_{v_{2}}(p, {\tilde s})|^{2}dp\leq 
$$
$$
(2\pi)^{2} e^{2aT}\|\widehat{G}(p)\|_{L^{\infty}({\mathbb R}^{2})}^{2}\int_{{\mathbb R}^{2}}
|\widehat{f}_{v_{1}}(p, {\tilde s})-\widehat{f}_{v_{2}}(p, {\tilde s})|^{2}dp\leq 
$$
$$
e^{2aT}\|G(x)\|_{L^{1}({\mathbb R}^{2})}^{2}\|F(v_{1}(x, {\tilde s}), x)-F(v_{2}(x, {\tilde s}), x)\|_
{L^{2}({\mathbb R}^{2})}^{2}.
$$
Let us use condition (\ref{lip}). Thus,
\begin{equation}
\label{lipl2}  
\|F(v_{1}(x, {\tilde s}), x)-F(v_{2}(x, {\tilde s}), x)\|_{L^{2}({\mathbb R}^{2})}\leq l
\|v_{1}(x, {\tilde s})-v_{2}(x, {\tilde s})\|_{L^{2}({\mathbb R}^{2})},
\end{equation}
so that
$$
\|e^{(t-{\tilde s})\{-p_{1}^{2}-|p_{2}|^{2s}+ib.p+a\}}2\pi \widehat{G}(p)
[\widehat{f}_{v_{1}}(p, {\tilde s})-\widehat{f}_{v_{2}}(p, {\tilde s})]\|_{L^{2}({\mathbb R}^{2})}\leq
$$
\begin{equation}
\label{eatgvl2}  
e^{aT}l\|G(x)\|_{L^{1}({\mathbb R}^{2})}\|v_{1}(x, {\tilde s})-v_{2}(x, {\tilde s})\|_{L^{2}({\mathbb R}^{2})}.   
\end{equation}
Inequalities (\ref{int0tetsf12}) and (\ref{eatgvl2}) give us
$$
\|\widehat{u_{1}}(p,t)-\widehat{u_{2}}(p,t)\|_{L^{2}({\mathbb R}^{2})}\leq
$$
$$
e^{aT}l\|G(x)\|_{L^{1}({\mathbb R}^{2})}\int_{0}^{T}
\|v_{1}(x, {\tilde s})-v_{2}(x, {\tilde s})\|_{L^{2}({\mathbb R}^{2})}d {\tilde s}.
$$
According to the Schwarz inequality, we have
$$
\int_{0}^{T}\|v_{1}(x, {\tilde s})-v_{2}(x, {\tilde s})\|_{L^{2}({\mathbb R}^{2})}d {\tilde s}\leq
$$
\begin{equation}
\label{sch2}  
\sqrt{\int_{0}^{T}\|v_{1}(x, {\tilde s})-v_{2}(x, {\tilde s})\|_{L^{2}({\mathbb R}^{2})}^{2}d {\tilde s}}\sqrt{T}.
\end{equation}
Therefore,
$$  
\|\widehat{u_{1}}(p,t)-\widehat{u_{2}}(p,t)\|_{L^{2}({\mathbb R}^{2})}\leq
$$
\begin{equation}
\label{u1hu2hl2}  
e^{aT}l\sqrt{T}\|G(x)\|_{L^{1}({\mathbb R}^{2})}\|v_{1}(x,t)-v_{2}(x,t)\|_
{L^{2}({\mathbb R}^{2}\times [0, T])}.
\end{equation}
This yields
$$
\|u_{1}(x,t)-u_{2}(x,t)\|_{L^{2}({\mathbb R}^{2}\times [0, T])}^{2}=\int_{0}^{T}
\|\widehat{u_{1}}(p,t)-\widehat{u_{2}}(p,t)\|_{L^{2}({\mathbb R}^{2})}^{2}dt\leq
$$
\begin{equation}
\label{u1u2l2v1v2}
e^{2aT}l^{2}T^{2}\|G(x)\|_{L^{1}({\mathbb R}^{2})}^{2}
\|v_{1}(x,t)-v_{2}(x,t)\|_{L^{2}({\mathbb R}^{2}\times [0, T])}^{2}.
\end{equation}  
Recall (\ref{u1u2hint0t}). Hence,
$$
p^{2}[\widehat{u_{1}}(p,t)-\widehat{u_{2}}(p,t)]=
\int_{0}^{t}
e^{(t-{\tilde s})\{-p_{1}^{2}-|p_{2}|^{2s}+ib.p+a\}}2\pi p^{2}\widehat{G}(p)
[\widehat{f}_{v_{1}}(p, {\tilde s})-\widehat{f}_{v_{2}}(p, {\tilde s})]d {\tilde s}.
$$
Clearly, the estimate on the norm 
$$
\|p^{2}[\widehat{u_{1}}(p,t)-\widehat{u_{2}}(p,t)]\|_{L^{2}({\mathbb R}^{2})}\leq
$$
\begin{equation}
\label{int0tetsf122}  
\int_{0}^{t}\|
e^{(t-{\tilde s})\{-p_{1}^{2}-|p_{2}|^{2s}+ib.p+a\}}2\pi p^{2}\widehat{G}(p)
[\widehat{f}_{v_{1}}(p, {\tilde s})-\widehat{f}_{v_{2}}(p, {\tilde s})]\|_{L^{2}({\mathbb R}^{2})}d {\tilde s}
\end{equation}
holds.
By virtue of (\ref{fub1}), we obtain the upper bound
$$
\|e^{(t-{\tilde s})\{-p_{1}^{2}-|p_{2}|^{2s}+ib.p+a\}}2\pi p^{2}\widehat{G}(p)
[\widehat{f}_{v_{1}}(p, {\tilde s})-\widehat{f}_{v_{2}}(p, {\tilde s})]\|_{L^{2}({\mathbb R}^{2})}^{2}=
$$
$$
(2\pi)^{2} \int_{{\mathbb R}^{2}}e^{-2(t-{\tilde s})p_{1}^{2}}e^{-2(t-{\tilde s})|p_{2}|^{2s}}e^{2(t-{\tilde s})a}|p^{2}\widehat{G}(p)|^{2}
|\widehat{f}_{v_{1}}(p, {\tilde s})-\widehat{f}_{v_{2}}(p, {\tilde s})|^{2}dp\leq 
$$
$$
(2\pi)^{2} e^{2aT}\|p^{2}\widehat{G}(p)\|_{L^{\infty}({\mathbb R}^{2})}^{2}\int_{{\mathbb R}^{2}}
|\widehat{f}_{v_{1}}(p, {\tilde s)}-\widehat{f}_{v_{2}}(p, {\tilde s})|^{2}dp\leq 
$$
$$
e^{2aT}\|\Delta_{x}G(x)\|_{L^{1}({\mathbb R}^{2})}^{2}
\|F(v_{1}(x, {\tilde s}), x)-F(v_{2}(x, {\tilde s}), x)\|_{L^{2}({\mathbb R}^{2})}^{2}.
$$
Let us recall formula (\ref{lipl2}). Thus,
$$
\|e^{(t-{\tilde s})\{-p_{1}^{2}-|p_{2}|^{2s}+ib.p+a\}}2\pi p^{2}\widehat{G}(p)
[\widehat{f}_{v_{1}}(p, {\tilde s})-\widehat{f}_{v_{2}}(p, {\tilde s})]\|_{L^{2}({\mathbb R}^{2})}\leq 
$$
\begin{equation}
\label{etsp2gpeat}
e^{aT}\|\Delta_{x} G(x)\|_{L^{1}({\mathbb R}^{2})}l\|v_{1}(x, {\tilde s})-v_{2}(x, {\tilde s})\|_
{L^{2}({\mathbb R}^{2})}.
\end{equation}
By means of (\ref{int0tetsf122}) along with (\ref{sch2}) and (\ref{etsp2gpeat}), we
arrive at
$$
\|p^{2}[\widehat{u_{1}}(p,t)-\widehat{u_{2}}(p,t)]\|_{L^{2}({\mathbb R}^{2})}\leq
$$
\begin{equation}
\label{p2u1hu2hpt}  
e^{aT}\sqrt{T}l\|\Delta_{x} G(x)\|_{L^{1}({\mathbb R}^{2})}
\|v_{1}(x,t)-v_{2}(x,t)\|_{L^{2}({\mathbb R}^{2}\times [0, T])}.
\end{equation}
Therefore,
$$
\|\Delta_{x} [u_{1}(x,t)-u_{2}(x,t)]\|_
{L^{2}({\mathbb R}^{2}\times [0, T])}^{2}=\int_{0}^{T}
\|p^{2}[\widehat{u_{1}}(p,t)-\widehat{u_{2}}(p,t)]\|_{L^{2}({\mathbb R}^{2})}^{2}dt\leq     $$
\begin{equation}
\label{u1u2l2v1v22}
e^{2aT}l^{2}T^{2}\|\Delta_{x} G(x) \|_{L^{1}({\mathbb R}^{2})}^{2}
\|v_{1}(x,t)-v_{2}(x,t)\|_{L^{2}({\mathbb R}^{2}\times [0, T])}^{2}.
\end{equation}  
From (\ref{u1u2hint0t}) we easily derive that
$$
\frac{\partial}{\partial t}[\widehat{u_{1}}(p,t)-\widehat{u_{2}}(p,t)]=
$$
$$
\{-p_{1}^{2}-|p_{2}|^{2s}+ib.p+a\}[\widehat{u_{1}}(p,t)-\widehat{u_{2}}(p,t)]+
2\pi \widehat{G}(p)[\widehat{f}_{v_{1}}(p,t)-\widehat{f}_{v_{2}}(p,t)].
$$
Then
$$
\Big\|\frac{\partial}{\partial t}[\widehat{u_{1}}(p,t)-\widehat{u_{2}}(p,t)]
\Big\|_{L^{2}({\mathbb R}^{2})}\leq a\|\widehat{u_{1}}(p,t)-\widehat{u_{2}}(p,t)\|_
{L^{2}({\mathbb R}^{2})}+
$$
$$      
\|b.p[\widehat{u_{1}}(p,t)-\widehat{u_{2}}(p,t)]\|_
{L^{2}({\mathbb R}^{2})}+
\||p_{2}|^{2s}[\widehat{u_{1}}(p,t)-\widehat{u_{2}}(p,t)]\|_{L^{2}({\mathbb R}^{2})}+    
$$
\begin{equation}
\label{ddtu1hptu2hpt}
\|p_{1}^{2}[\widehat{u_{1}}(p,t)-\widehat{u_{2}}(p,t)]\|_{L^{2}({\mathbb R}^{2})}+ 
2\pi \|\widehat{G}(p)[\widehat{f}_{v_{1}}(p,t)-\widehat{f}_{v_{2}}(p,t)]\|_
{L^{2}({\mathbb R}^{2})}.
\end{equation}  
Recall (\ref{u1hu2hl2}). Hence, the first term in the right side of
(\ref{ddtu1hptu2hpt}) can be estimated from above by
\begin{equation}
\label{au1hu2hl2}  
a \sigma e^{aT}\sqrt{T}l\|v_{1}(x,t)-v_{2}(x,t)\|_{L^{2}({\mathbb R}^{2}\times [0, T])}
\end{equation}
with $\sigma$ introduced in (\ref{q}).
We obtain the upper bound on the norm as
$$
\|b.p[\widehat{u_{1}}(p,t)-\widehat{u_{2}}(p,t)]\|_{L^{2}({\mathbb R}^{2})}^{2}\leq
$$
$$
|b|_{{\mathbb R}^{2}}^{2}\Bigg[\int_{|p|\leq 1}p^{2}|\widehat{u_{1}}(p,t)-\widehat{u_{2}}(p,t)|^{2}dp+
\int_{|p|>1}p^{2}|\widehat{u_{1}}(p,t)-\widehat{u_{2}}(p,t)|^{2}dp \Bigg]\leq 
$$
$$  
|b|_{{\mathbb R}^{2}}^{2}[\|\widehat{u_{1}}(p,t)-\widehat{u_{2}}(p,t)\|_{L^{2}({\mathbb R}^{2})}^{2}+
\|p^{2}[\widehat{u_{1}}(p,t)-\widehat{u_{2}}(p,t)]\|_{L^{2}({\mathbb R}^{2})}^{2}].
$$
Let us use inequalities (\ref{u1hu2hl2}) and  (\ref{p2u1hu2hpt}). They
imply the estimate from above on the second term in the right side of
(\ref{ddtu1hptu2hpt}) equal to
\begin{equation}
\label{bu1hu2hl2}  
|b|_{{\mathbb R}^{2}}\sigma e^{aT}\sqrt{T}l\|v_{1}(x,t)-v_{2}(x,t)\|_{L^{2}({\mathbb R}^{2}\times [0, T])}.
\end{equation}
Evidently,
$$
\||p_{2}|^{2s}[\widehat{u_{1}}(p,t)-\widehat{u_{2}}(p,t)]\|_{L^{2}({\mathbb R}^{2})}^{2}\leq
$$
$$
\int_{|p|\leq 1}|p|^{4s}|\widehat{u_{1}}(p,t)-\widehat{u_{2}}(p,t)|^{2}dp+
\int_{|p|>1}|p|^{4s}|\widehat{u_{1}}(p,t)-\widehat{u_{2}}(p,t)|^{2}dp\leq 
$$
$$  
\|\widehat{u_{1}}(p,t)-\widehat{u_{2}}(p,t)\|_{L^{2}({\mathbb R}^{2})}^{2}+
\|p^{2}[\widehat{u_{1}}(p,t)-\widehat{u_{2}}(p,t)]\|_{L^{2}({\mathbb R}^{2})}^{2}.
$$
By means of (\ref{u1hu2hl2}) and  (\ref{p2u1hu2hpt}), the third term in the
right side of (\ref{ddtu1hptu2hpt}) can be bounded from above by
\begin{equation}
\label{bu1hu2hl3}  
\sigma e^{aT}\sqrt{T}l\|v_{1}(x,t)-v_{2}(x,t)\|_{L^{2}({\mathbb R}^{2}\times [0, T])}.
\end{equation}
Similarly,
$$
\|p_{1}^{2}[\widehat{u_{1}}(p,t)-\widehat{u_{2}}(p,t)]\|_{L^{2}({\mathbb R}^{2})}^{2}\leq
$$
$$
\int_{|p|\leq 1}|p|^{4}|\widehat{u_{1}}(p,t)-\widehat{u_{2}}(p,t)|^{2}dp+
\int_{|p|>1}|p|^{4}|\widehat{u_{1}}(p,t)-\widehat{u_{2}}(p,t)|^{2}dp\leq 
$$
$$  
\|\widehat{u_{1}}(p,t)-\widehat{u_{2}}(p,t)\|_{L^{2}({\mathbb R}^{2})}^{2}+
\|p^{2}[\widehat{u_{1}}(p,t)-\widehat{u_{2}}(p,t)]\|_{L^{2}({\mathbb R}^{2})}^{2}.
$$
Therefore, the fourth term in the right side of (\ref{ddtu1hptu2hpt}) can be estimated
from above by expression (\ref{bu1hu2hl3}) as well.
Inequalities (\ref{fub}) and (\ref{lipl2}) enable us to derive that
$$
(2\pi)^{2}\int_{{\mathbb R}^{2}}|\widehat{G}(p)|^{2}
|\widehat{f}_{v_{1}}(p,t)-\widehat{f}_{v_{2}}(p,t)|^{2}dp\leq
$$
$$
(2\pi)^{2}
\|\widehat{G}(p)\|_{L^{\infty}({\mathbb R}^{2})}^{2}\int_{{\mathbb R}^{2}}
|\widehat{f}_{v_{1}}(p,t)-\widehat{f}_{v_{2}}(p,t)|^{2}dp\leq
$$
$$
\|G(x)\|_{L^{1}({\mathbb R}^{2})}^{2}\|F(v_{1}(x,t), x)-F(v_{2}(x,t), x)\|_
{L^{2}({\mathbb R}^{2})}^{2}\leq
$$
$$
\|G(x)\|_{L^{1}({\mathbb R}^{2})}^{2}l^{2}\|v_{1}(x,t)-v_{2}(x,t)\|_{L^{2}({\mathbb R}^{2})}^{2}.  
$$
Thus, the fifth term in the right side of (\ref{ddtu1hptu2hpt}) can be bounded
from above by
\begin{equation}
\label{bu1hu2hl4}  
\sigma l\|v_{1}(x,t)-v_{2}(x,t)\|_{L^{2}({\mathbb R}^{2})}.
\end{equation}
Let us combine (\ref{au1hu2hl2}), (\ref{bu1hu2hl2}), (\ref{bu1hu2hl3}) and
(\ref{bu1hu2hl4}). This yields
$$
\Big\|\frac{\partial}{\partial t}[\widehat{u_{1}}(p,t)-\widehat{u_{2}}(p,t)]
\Big\|_{L^{2}({\mathbb R}^{2})}\leq
$$
$$
\sigma e^{aT}\sqrt{T}l\{a+|b|_{{\mathbb R}^{2}}+2\}\|v_{1}(x,t)-v_{2}(x,t)\|_
{L^{2}({\mathbb R}^{2}\times [0, T])}+
\sigma l\|v_{1}(x,t)-v_{2}(x,t)\|_{L^{2}({\mathbb R}^{2})}.
$$
This allows us to estimate the norm as
$$
\Big\|\frac{\partial}{\partial t}(u_{1}(x,t)-u_{2}(x,t))\Big\|_
{L^{2}({\mathbb R}^{2}\times [0, T])}^{2}=\int_{0}^{T}
\Big\|\frac{\partial}{\partial t}[\widehat{u_{1}}(p,t)-\widehat{u_{2}}(p,t)]
\Big\|_{L^{2}({\mathbb R}^{2})}^{2}dt\leq           
$$
\begin{equation}
\label{ddtu1u2l2}  
2{\sigma}^{2}l^{2}[e^{2aT}T^{2}\{a+|b|_{{\mathbb R}^{2}}+2\}^{2}+1]
\|v_{1}(x,t)-v_{2}(x,t)\|_{L^{2}({\mathbb R}^{2}\times [0, T])}^{2}.
\end{equation}
Let us introduce a finite, positive constant 
\begin{equation}
\label{alf}
\alpha_{a, b, T}:=e^{2aT}T^{2}\{a+|b|_{{\mathbb R}^{2}}+2\}^{2}+1.
\end{equation}
Recall the definition of the norm (\ref{122n}). We combine upper
bounds (\ref{u1u2l2v1v2}), (\ref{u1u2l2v1v22}) and (\ref{ddtu1u2l2}). 
A trivial calculation gives us that
\begin{equation}
\label{contr}
\|u_{1}-u_{2}\|_{W^{(1, 2), 2}({\mathbb R}^{2}\times [0, T])}\leq
\sigma l\sqrt{T^{2}e^{2aT}+2\alpha_{a, b, T}}
\|v_{1}-v_{2}\|_{W^{(1, 2), 2}({\mathbb R}^{2}\times [0, T])}.
\end{equation}
Note that the constant in the right side of (\ref{contr}) is less than one
due to condition (\ref{qlt}). This means that under the stated assumptions equation
(\ref{aux}) defines the map
$$
T_{a, b}: W^{(1, 2), 2}({\mathbb R}^{2}\times [0, T])\to W^{(1, 2), 2}
({\mathbb R}^{2}\times [0, T]),
$$
which is a strict contraction. Its unique fixed
point $w(x,t)$ is the only solution of
problem (\ref{h}), (\ref{ic}) in
$W^{(1, 2), 2}({\mathbb R}^{2}\times [0, T])$. \hfill\lanbox

\bigskip

\noindent
{\it Proof of Corollary 1.4.} The statement of the Corollary is valid
because the constant in the right side of bound
(\ref{contr}) is independent of the initial condition (\ref{ic}) (see e.g. ~\cite{EZ01}). Hence,
problem (\ref{h}), (\ref{ic}) admits a unique solution
$w(x,t)\in W^{(1, 2), 2}({\mathbb R}^{2}\times {\mathbb R}^{+})$. Suppose
that $w(x,t)$ is trivial for $x\in {\mathbb R}^{2}$ and $t\in {\mathbb R}^{+}$.
Then we arrive at the contradiction to the given condition that
$\hbox{supp}\widehat{F(0, x)}\cap \hbox{supp}\widehat{G}$ is a set of
nonzero Lebesgue measure in ${\mathbb R}^{2}$. \hfill\lanbox

\bigskip


\centerline{\bf 3. Acknowledgement}

\bigskip

\noindent
V.V. is grateful to Israel Michael Sigal for the partial support by the
NSERC grant NA 7901.

\end{document}